[1994/12/01]
\documentclass[a4paper,12pt]{article}
\usepackage{amsfonts,amssymb,amsmath,amsthm}
\usepackage{eurosym}

\input{tcilatex}
\chardef\bslash=`\\ 

\newtheorem{st}{Statement}[section]

\newtheorem{propo}[st]{Proposition}
\newtheorem{cor}[st]{Corollary}
\newtheorem{thm}[st]{Theorem}
\newtheorem{lemm}[st]{Lemma}

\newtheorem{fact}[st]{Fact}

\theoremstyle{definition}

\newtheorem{cl}[st]{Claim}
\newtheorem{rem}[st]{Remark}
\newtheorem{quest}[st]{Question}

\theoremstyle{remark}

\newcommand{\eval}[2][\right]{\relax
  \ifx#1\right\relax \left.\fi#2#1\rvert}

\def\Level{{\rm Level\,}}

\def\cof{{\rm cof\,}}

\def\Level{{\rm Level\,}}
\def\lim{{\rm lim\,}}
\def\fin{{\rm fin\,}}

\def\Proof:{ \vspace{-1.5mm} {\noindent\it Proof:}}

\begin{document}

\title{\vspace{0mm} \vspace{0mm} The Rudin--Keisler ordering of P-points
under $\mathfrak{b} = \mathfrak{c}$}
\author{Andrzej Starosolski}
\date{\today}
\maketitle

%

\begin{abstract}
M. E. Rudin (1971) proved, under CH, that for each P-point $p$ there exists a
P-point $q$ strictly RK-greater than $p$. This result was proved under ${\mathfrak{p}=%
\mathfrak{c}}$ by A. Blass (1973), who also showed that each RK-increasing $%
\hspace{0.5mm}\omega $-sequence \hspace{0.4mm} of \hspace{0.4mm} P-points%
\hspace{0.4mm} is \hspace{0.4mm} upper \hspace{0.4mm}bounded \hspace{0.4mm}
by a P-point, and that there is an order embedding of the real line into the
class of P-points with respect to the RK-ordering. In this paper, the
results cited above are proved under the (weaker) assumption that $\mathfrak{%
b}=\mathfrak{c}$.

A. Blass asked in 1973 which ordinals can be embedded in the set of
P-points, and pointed out that such an ordinal cannot be greater than $%
\mathfrak{c}^{+}$. In this paper it is proved, under $\mathfrak{b}=\mathfrak{c}$, that for each
ordinal $\alpha < \mathfrak{c}^{+}$, there is an order embedding of $%
\alpha$ into P-points. It is also proved, under $\mathfrak{b}=\mathfrak{c}$,
that there is an embedding of the long line into P-points.
%
\end{abstract}

\section{Introduction}


In \cite{MERudin}, M. E. Rudin proved that, under CH, for each P-point $p$ there
exists a P-point $q$ strictly RK-greater than $p$. A. Blass showed the same \cite{Blass1}
assuming that $\mathfrak{p}=\mathfrak{c}$  \footnote{
Actually, all results from \cite{Blass1} quoted in this paper were stated
under MA, but the proofs also work under $\mathfrak{p}=\mathfrak{c}$, as
pointed out by A. Blass in \cite{Blass2};
definitions of $\mathfrak{p}, \mathfrak{c}$ and few others cardinal invariants are recalled on page 3.};
moreover, he proved that if $\mathfrak{p}=\mathfrak{c}$, then each RK-increasing $\omega$-sequence of P-points is upper bounded by a P-point, and there exists an order-embedding of the real line into the class of P-points with respect to the RK-ordering.
Since then, the RK-ordering of P-points has been thoroughly investigated,
however, most of the obtained results were proved under $MA_{\sigma -\mathrm{%
centr}}$, or stronger assumptions  \footnote{ Note that
A. Blass asks \cite[Question 5]{Blass1}: What can be proved about P-points
without using MA?}, usually with complicated proofs and using sophisticated
techniques. We prove the results mentioned above under $\mathfrak{b}=%
\mathfrak{c}$. Perhaps more importantly, we present a method of proof that
turns out be effective in the study of P-points under $\mathfrak{b}=%
\mathfrak{c}$. The ideas used in the present paper were originally presented
in an unpublished paper \cite{Star_RK_P-hier}, where the RK-ordering
concerned the ultrafilters in the classes of the so-called P-hierarchy, the
first class of which coincides with that of P-points. The method is based on
the use of contours and quasi-subbases, which enables us to employ
surprisingly concise arguments, in contrast with the approaches of some
other papers on similar topics.

After a scrutiny of mechanisms underlying our proofs, we introduce an
apparently new cardinal invariant $\mathfrak{q}$, the use of which enables
us to weaken the set-theoretic assumptions of most of our results. Finally,
we show that $\mathfrak{q}$ is an instance of a general method of
constructing useful variants of cardinal invariants.

In a recent paper, D. Raghavan and S. Shelah \cite{Raghavan-Shelah} proved
(under $\mathfrak{p}=\mathfrak{c}$) that there is an order-embedding of ${%
\mathcal{P}}(\omega )/\fin$ into the set of P-points ordered by $\geq _{RK}$,
and gave a short review of earlier results concerning embeddings of
different orders into the class of P-points.

A. Blass also asked \cite[Question 4]{Blass1} which ordinals can be embedded
in the set of P-points, and pointed out that such an ordinal cannot be
greater than $\mathfrak{c}^{+}$. We show that under $\mathfrak{b}=\mathfrak{c}$
each ordinal less than $\mathfrak{c}^{+}$ is order-embeddable into
P-points. A recent paper by B. Kuzeljevic and D. Raghavan \cite{Kuz-Rag}
answers the question of the embedding of ordinals into P-points under MA.

\section{Tools}

A free ultrafilter $u$ is a \emph{P-point} if and only if, for each
partition $(V_{n})_{n<\omega }$ of $\omega $, there exists a set $U\in u$
such that either $U\subset V_{n}$ for some $n<\omega $ or else $U\cap V_{n}$
is finite for all $n<\omega $. A filter $\mathcal{F}$ is said to be {\it Rudin-Keisler greater} (\emph{%
RK-greater}) than a filter $\mathcal{G}$ (written as $\mathcal{F}$ $\geq _{RK}%
\mathcal{G}$) if there exists a map $h$ such that
$G \in \mathcal{G}$ if and only if $h^{-1}(G)\in \mathcal{F}$. Let
\begin{equation}
{\mathcal{W}}=\{W_{n}:{n<\omega }\}  \label{part}
\end{equation}%
be a \emph{partition} of a subset of $\omega $ into infinite sets. A filter $%
\mathcal{K}$ is called a \emph{contour} if there exists a partition $%
\mathcal{W}$ such that $W\in \mathcal{K}$ if and only if there is a cofinite set $%
I\subset \omega $ such that $W\cap W_{n}$ is cofinite on $W_{n}$ for each $%
n\in I$. We call $\mathcal{K}$ a contour of $\mathcal{W}$, and denote $%
\mathcal{K}=\int \mathcal{W}$. \footnote{See \cite%
{Dol-Myn1} and the last section of \cite{DSW1} for a systematic presentation
of contours.}

A fundamental property used in the present paper is the following
reformulation of \cite[Proposition 2.1]{Star-P-hier}.

\begin{propo}
\label{A0} 
A free ultrafilter is a P-point if and only if it does not include a contour. %
\end{propo}

As usual, $\mathfrak{c}$ denotes the cardinality of the continuum.
If $f,g\in {^{\omega }\omega }$, then we say that $f$ \textit{%
dominates} $g$ (and write $f\geq ^{\ast }g$) if $f(n)\leq g(n)$ for almost
all $n<\omega $. We say that a family $\mathcal{F}$ of ${^{\omega }\omega }$
functions is \textit{unbounded} if there is no $g\in {^{\omega }\omega }$
that dominates all functions $f\in \mathcal{F}$. The minimal cardinality of
an unbounded family is the \textit{bounding number} $\mathfrak{b}$. We also
say that a family $\mathcal{F}\subset {^{\omega }\omega }$ is \textit{%
dominating} if, for each $g\in {^{\omega }\omega }$, there is some $f\in
\mathcal{F}$ that dominates $g$. The \textit{dominating number} $\mathfrak{d}
$ is the minimal cardinality of a dominating family. The \textit{%
pseudointersection number} $\mathfrak{p}$ is the minimal cardinality of a
free filter without \textit{pseudointersection}, which is a set almost included in each
element of the filter. Finally, the \textit{ultrafilter number} $\mathfrak{u%
}$ is the minimal cardinality of a base of a free ultrafilter. It is well
known that $\mathfrak{b}\leq \mathfrak{d}\leq \mathfrak{c}$, and $\mathfrak{p%
}\leq \mathfrak{b}\leq \mathfrak{u}\leq \mathfrak{c}$, and that there are
models for which $\mathfrak{p}<\mathfrak{b}$
(see, for example, \cite{vanDouwen}).

If $\mathcal{A}$ and $\mathcal{B}$ are families of sets, then we say that $%
\mathcal{A}$ and $\mathcal{B}$ \emph{mesh} (in symbols, $\mathcal{A}\#%
\mathcal{B}$) if $A\cap B\neq \varnothing $ for each $A\in \mathcal{A}$ and $%
B\in \mathcal{B}$, and we abridge $\left\{ A\right\} \#\mathcal{B}$ to $A\#%
\mathcal{B}$.

If a family $\mathcal{A}$ has the finite intersection property (fip), then we denote
by $\langle \mathcal{A}\rangle $ the filter generated by $\mathcal{A}$. A
family $\mathcal{A}$ of sets is said to be \emph{free} if there exists a
free filter $\mathcal{F}$ such that $\mathcal{A}\subset \mathcal{F}$. \footnote{A filter $\mathcal{F}$ is called \emph{free}
if $\bigcap\nolimits_{F\in \mathcal{F}}F=\varnothing .$} Consequently, each
free family has the finite intersection property. Let $A$ be an infinite
subset of $\omega $. A filter $\mathcal{F}$, on $\omega $, is said to be
\emph{cofinite on}\textit{\ }$A$ whenever a set $U\in \mathcal{F}$ if and only if $%
A\setminus U$ is finite. A filter $\mathcal{F}$ is said to be \emph{cofinite} if\emph{\ }it is
cofinite on some $A$. The class of cofinite filters is denoted by $\mathbb{C}%
\mathrm{of}$. It is well-known that a filter is free on $X$ if and only if
it includes the cofinite filter of $X$.

In the sequel we will often need to know whether a family $\mathcal{A}$ with
the finite intersection property is included in a free ultrafilter, which is
equivalent to the property that the filter $\langle \mathcal{A}\rangle $
admits a finer free ultrafilter. Here is a simple criterion, because of
which, free families are also known as families with the \textit{strong
finite intersection property}.

\begin{fact}[Folklore]
A filter $\mathcal{F}$ is included in a free ultrafilter if and only if all elements of
the filter $\mathcal{F}$ are infinite.
\end{fact}

\begin{proof}
Let $\mathcal{F}$ be a filter included in an ultrafilter $\mathcal{U}$. The
ultrafilter $\mathcal{U}$ is free if and only if it includes the cofinite
filter $\mathcal{C}$. So $\mathcal{F}\subset \mathcal{U}$ and $\mathcal{C}%
\subset \mathcal{U},$ hence $\mathcal{F}$\#$\mathcal{C}$. In other words, a
filter $\mathcal{F}$ is included in a free ultrafilter whenever $\mathcal{F}$%
\#$\mathcal{C}$, equivalently whenever all the elements of $\mathcal{F}$ are
infinite.
\end{proof}

We say that families $\mathcal{A}$ and $\mathcal{B}$ \emph{interact} (in symbols, $%
\mathcal{A}\ddagger \mathcal{B}$) if $\mathcal{A}\cup \mathcal{B}$ is free.
It is clear that $\mathcal{A}\ddagger \mathcal{B}$ implies $\mathcal{A}\#%
\mathcal{B}$. 
If $\mathcal{A}$ and $\mathcal{B}$ interact, then $A\cap B$ is infinite for
every $A\in \mathcal{A}$ and $B\in \mathcal{B}$. There exist free families $%
\mathcal{A}$ and $\mathcal{B}$ with infinite $A\cap B$ for every $A\in
\mathcal{A}$ and $B\in \mathcal{B}$, which do not interact.
\footnote{For example, let $X_{0},X_{1},\ldots ,X_{6}$ be
disjoint infinite sets. Let $A_{\alpha ,\ldots ,\gamma }=X_{\alpha }\cup
\ldots \cup X_{\gamma }.$ Consider $\mathcal{A}=\left\{
A_{1,2,3},A_{3,4,5}\right\} ,\mathcal{B}=\left\{ A_{1,4,6},A_{2,5,6}\right\}
.$}

\textbf{A relation between sets and functions.} Let $\mathcal{W}$ be a
partition \nolinebreak (\ref{part}). For each $n<\omega $, let $\left( w_{k}^{n}\right)
_{k<\omega }$ be an increasing sequence such that%
\begin{equation*}
W_{n}=\left\{ w_{k}^{n}k<\omega \right\} .
\end{equation*}%
For each $f\in {^{\omega }\omega }$ and $m<\omega $, let%
\begin{equation}
\limfunc{E}\nolimits_{\mathcal{W}}\left( f,m\right) =\left\{
w_{k}^{n}f\left( n\right) \leq k,m\leq n\right\}  \label{epi}
\end{equation}

If $F\in \int \mathcal{W}$, then, by definition, there exists the least $%
n_{F}<\omega $ such that $W_{n}\setminus F$ is finite for each $n\geq n_{F}$%
. Now, for each $n\geq n_{F}$, there exists a minimal $k_{n}<\omega $ such
that $w_{k}^{n}\in F$ for each $k\geq k_{n}$. Let ${f\hspace{-1.8mm}f}_{F}$
denote the set of those functions $f$ for which
\begin{equation}
n\geq n_{F}\Longrightarrow f\left( n\right) =k_{n}.  \label{such}
\end{equation}%
Then $\limfunc{E}\nolimits_{\mathcal{W}}\left( f,n_{F}\right) $ is the same
for each $f\in {f\hspace{-1.8mm}f}_{F}$, it is the largest set of the form
(\ref{epi}) included in $F$. Sure enough, $\limfunc{E}\nolimits_{\mathcal{W}%
}\left( f_{F},n_{F}\right) \in \int \mathcal{W}$.

Conversely, for every function $f\in {^{\omega }\omega }$, we define a
family $\mathcal{W}_{f}$ of subsets of $\omega $ as follows: $F\in \mathcal{W%
}_{f}$ if there is $n_{F}<\omega $ such that $F=\limfunc{E}\nolimits_{%
\mathcal{W}}\left( f_{F},n_{F}\right) $. Therefore, we can state the
following.

\begin{propo}
The family $\bigcup\nolimits_{f\in {^{\omega }\omega }}\mathcal{W}_{f}$ is a
base of $\int \mathcal{W}$.$\bigskip $
\end{propo}

\textbf{Quasi-subbases.\textit{\ }}We say that a family $\mathcal{A}$ is \emph{finer}
than $\mathcal{B}$ if $\langle \mathcal{B}\rangle \subset \langle \mathcal{A}%
\rangle $. Moreover, $\mathcal{A}$ is called a \emph{quasi-subbase} of (a
filter) $\mathcal{F}$ if there exists a countable family $\mathcal{C}$ such
that $\langle \mathcal{A}\cup \mathcal{C}\rangle =\mathcal{F}$. Accordingly,
$\mathcal{A}$ is \emph{quasi-finer} than $\mathcal{B}$ if there exists a
countable family $\mathcal{C}$ such that $\langle \mathcal{B}\rangle \subset
\langle \mathcal{A}\cup \mathcal{C}\rangle $.

%
%
\vspace{3mm} If $\mathcal{W}$ is a partition (\ref{part}), then for each $%
i<\omega $, let
\begin{equation*}
{\widetilde{W}_{i}=}\bigcup_{n\geq i}W_{n}\text{ and }\widetilde{\mathcal{W}}%
=\{{\widetilde{W}_{i}}:i<\omega \}.
\end{equation*}

\begin{propo}
\label{AAA} Let $\mathcal{W}$ be a partition and let $\mathcal{A}$ be a free
family. Then the following are equivalent:

(1) $\mathcal{A}$ is quasi-finer than $\int \mathcal{W}$,

(2) there exists a set $D$ such that $\mathcal{A}\cup \widetilde{\mathcal{W}}%
\cup \{D\}$ is free, and $\int {\mathcal{W}}\subset \langle \mathcal{A}\cup
\widetilde{\mathcal{W}}\cup \{D\}\rangle $.
\end{propo}

\begin{proof}
The implication $2\Rightarrow 1$ is evident. We will show $1\Rightarrow 2$. Suppose
the contrary, 
and let $\mathcal{B}$ be a countable family
of sets such that $\int \mathcal{W}\subset \langle \mathcal{A}\cup \mathcal{B%
}\rangle $. Taking finite intersections $\bigcap_{i\leq n}B_{i}$ instead of $%
B_{n},$ we obtain a decreasing sequence so that, without loss of generality,
we can assume that $\mathcal{B}=\{B_{n}\}_{n<\omega }$ is decreasing.
Since (2) is false, for each $n$ there exists $A_{n}\in \int \mathcal{W}$ such that $%
A_{n}\not\in \langle \mathcal{A}\cup \widetilde{\mathcal{W}}\cup
\{B_{n}\}\rangle $.

Without loss of generality, for each $n$ there is $k(n)\geq n$ such that $%
A_{n}\cap W_{i}$ is empty for all $i<k(n)$ and $W_{i}\setminus A_{n}$ is
finite for all $i\geq k(n)$. Define $A_{\infty }=\bigcup_{i<\omega }\left(
\bigcap_{\{n:k(n)\leq i\}}A_{n}\cap W_{i}\right) $ and note that $A_{\infty
}\in \int \mathcal{W}$.

We will show that $A_{\infty }\not\in \langle \mathcal{A}\cup \mathcal{B}%
\cup \widetilde{\mathcal{W}}\rangle \supset \langle \mathcal{A}\cup \mathcal{%
B}\rangle $. For this purpose, it suffices to show that $A_{\infty }\not\in
\langle \mathcal{A}\cup \widetilde{\mathcal{W}}\cup \{B_{n}\}\rangle $ for
each $n<\omega $. Indeed, note that $A_{\infty }\subset A_{n}\cup \widetilde{%
W}_{k(n)}^{c}$ for each $n<\omega $. From $(A_{n}\cup \widetilde{W}%
_{k(n)}^{c})\cap \widetilde{W}_{k(n)}\subset A_{n}\not\in \langle \mathcal{A}%
\cup \widetilde{\mathcal{W}}\cup \{B_{n}\}\rangle $, we infer that $%
(A_{n}\cup \widetilde{W}_{k(n)}^{c})\not\in \langle \mathcal{A}\cup
\widetilde{\mathcal{W}}\cup \{B_{n}\}\rangle $, and so $A_{\infty }\not\in
\langle \mathcal{A}\cup \widetilde{\mathcal{W}}\cup \{B_{n}\}\rangle$
\end{proof}

\begin{rem}
\label{AAA111} Let $\mathcal{A}$ be finer than $\int \mathcal{W}$. Then
there exists a partition $\mathcal{V}$ such that $\int \mathcal{V}\subset
\langle \mathcal{A}\cup \widetilde{\mathcal{V}}\rangle $.
\end{rem}

\begin{proof}
Take $V_i$ to be the $i$-th infinite $W_j \cap D$, with $D$ from Proposition %
\ref{AAA}(2)).
\end{proof}

\begin{thm}
\label{DDD} Let $(\mathcal{A}_{\alpha })_{\alpha <\beta <\mathfrak{b}}$ be
an increasing sequence of families of sets such that $\mathcal{A}_{\alpha }$
is not quasi-finer than a contour. Then, $\bigcup_{\alpha <\beta }\mathcal{A}%
_{\alpha }$ is not quasi-finer than a contour.
\end{thm}

\begin{proof}
By way of contradiction, $\{\mathcal{A}_{\alpha }\}_{\alpha <\beta <%
\mathfrak{b}}$ be as in the assumptions, and let $\int \mathcal{W}$ be a
contour and $\mathcal{B}$ be a countable family of sets such that $\int
\mathcal{W}\subset \langle \mathcal{A}\cup \mathcal{B}\rangle $. By
Proposition \ref{AAA} and Remark \ref{AAA111}, without loss of generality,
we may assume that $\mathcal{B}=\{B_{n}\}_{n<\omega }=\widetilde{\mathcal{W}}
$ is decreasing. Let $\mathcal{C}_{\alpha }=\langle \mathcal{A_{\alpha }}%
\cup \mathcal{B}\rangle $. Clearly $\mathcal{C}_{\alpha }$ does not include
$\int \mathcal{W}$ for each $\alpha <\beta $. Thus, for each $\alpha <\beta $%
, there exists a set $D_{\alpha }\in \int \mathcal{W}$ such that $D_{\alpha
}\not\in \mathcal{C}_{\alpha }$. Let $g_{\alpha }\in {f\hspace{-1.8mm}f}%
_{D_{\alpha }}$ for each $\alpha <\beta $. As $\beta <\mathfrak{b}$, the
family $\{g_{\alpha }\}_{\alpha <\beta }$ is bounded by some function $g$.
Let $G\in \mathcal{W}_{g}$. We will show that $G\not\in \bigcup_{\alpha
<\beta }\mathcal{C}_{\alpha }$ hence $G\not\in \mathcal{C}_{\alpha }$ for
each $\alpha <\beta $. Suppose not, and let $\alpha _{0}$ be a witness. By
construction, there exists $n_{0}<\omega $ such that $G\subset
B_{n_{0}}^{c}\cup D_{\alpha _{0}}$. As $(D_{\alpha _{0}}\cup
B_{n_{0}}^{c})\cap B_{n_{0}}\subset D_{\alpha _{0}}\not\in \mathcal{C}%
_{\alpha _{0}}$, it follows that $D_{\alpha _{0}}\cup B_{n_{0}}^{c}\not\in
\mathcal{C}_{\alpha _{0}}$, and so $G\not\in \mathcal{C}_{\alpha _{0}}$.
\end{proof}

\begin{cor}
\label{DDD111} ($\mathfrak{b}=\mathfrak{c}$) If $\mathcal{A}$ is a free
family of subsets of $\omega $ that is not quasi-finer than a contour, then
there exists a P-point $p$ such that $\mathcal{A}\subset p$.
\end{cor}

\begin{proof}
We arrange all contours in a sequence $(\int {\mathcal{W}}_{\alpha
})_{\alpha <\mathfrak{b}}$. We will build an increasing $\mathfrak{b}$%
-sequences $({\mathcal{F}}_{\alpha })_{\alpha <\mathfrak{b}}$ of filters
such that:

1) Each ${\mathcal{F}}_{\alpha }$ is generated by $\mathcal{A}$ together
with some family of sets of cardinality $<\mathfrak{b}$;

2) ${\mathcal{F}}_0 = \mathcal{A}$;

3) For each $\alpha <\mathfrak{b}$, there exists $F\in {\mathcal{F}}_{\alpha
+1}$ such that $F^{c}\in \int {\mathcal{W}}_{\alpha }$;

4) For every limit $\alpha $ let ${\mathcal{F}}_{\alpha }=\bigcup_{\beta
<\alpha }{\mathcal{F}}_{\beta } $.

The existence of such a sequence follows using standard induction by Theorem %
\ref{DDD}. Now it suffices to take any ultrafilter extending $%
\bigcup_{\alpha <\mathfrak{b}}{\mathcal{F}}_{\alpha }$ and note that, by
Proposition \ref{A0}, it is a P-point.
\end{proof}

Let us recall a well-known theorem (see, for example \cite[Corollary 1]%
{Blass1}).

\begin{thm}
\label{MMM111} Let $u$ be an ultrafilter. If $f(u)=_{RK}u$, then there
exists $U\in u$ such that $f$ is one-to-one on $U$.
\end{thm}

\section{Applications: RK-ordering of P-points}

M. E. Rudin \cite{MERudin} proved that, under CH, for each P-point $p$ there
exists a P-point $q$ strictly RK-greater than $p$. Some years later, A. Blass \cite[%
Theorem 6]{Blass1} proved this theorem under $\mathfrak{p}=\mathfrak{c}$.

\begin{thm}
\label{III} ($\mathfrak{b}=\mathfrak{c}$) \footnote{
Note that all theorems in this section are proved, in fact, under (possibly)
weaker assumptions, what we will discuss in detail in Section 4.} If $p$ is
a P-point, then there exists a P-point $q$ that is strictly Rudin--Keisler
greater than $p$.
\end{thm}

\begin{proof}
Let $f\in {^{\omega }\omega }$ be a finite-to-one function such that
\begin{equation*}
\limsup\nolimits_{n\in P}\mathrm{card\,}(f^{-1}(n))=\infty
\end{equation*}
holds for all $P\in p$. We define a family $\mathcal{A}$ as follows: $A\in
\mathcal{A}$ if and only if there exist $i<\omega $ and $P\in p$ such that $\mathrm{%
card\,}(f^{-1}(n)\setminus A)<i$ for each $n\in P$. Then, Theorem \ref%
{MMM111} ensures that the ultrafilters we are building are strictly
RK-greater than $p $.

We claim that $\left\{ f^{-1}(p)\right\} \cup \mathcal{A}$ is not
quasi-finer than any contour. Suppose not, and take a witness $\int \mathcal{%
W}$. From Remark \ref{AAA111}, without loss of generality, we may assume
that $\int {\mathcal{W}}\subset \langle f^{-1}(p)\cup \mathcal{A}\cup
\widetilde{\mathcal{W}}\rangle $. Consider two cases:

Case 1: There exists a sequence $(B_{n})_{n<\omega }$ and a strictly
increasing $k\in {^{\omega }\omega }$ such that $B_{n}\subset W_{k(n)}$, $%
f(B_{n})\not\in p$, and $B\#f^{-1}(p)\cup \mathcal{A} \cup \widetilde{
\mathcal{W}}$, where $B=\bigcup_{n<\omega }B_{n}$. Take a sequence $%
(f(\bigcup_{i\leq n}B_{n}))_{n<\omega }$. This is an increasing sequence,
and it is clear that $\bigcup_{n<\omega }f(\bigcup_{i\leq n}B_{n})=f(B)\in p$%
. Make a partition of $f(B)$ by taking $f(\bigcup_{n\leq i+1}B_{n})\setminus
f(\bigcup_{n\leq i}B_{n})$ for $i<\omega $. As $p$ is a P-point, there
exists some $P\in p$ such that $P\cap (f(\bigcup_{n\leq i+1}B_{n})\setminus
f(\bigcup_{n\leq i}B_{n}))$ is finite for all $i<\omega $, and thus $%
f^{-1}(P)\cap B_{n}$ is finite for all $i<\omega $. Therefore, $%
f^{-1}(P)\cap W_{i}\cap B$ is finite, and thus $(f^{-1}(P)\cap B)^{c}\in
\int \mathcal{W}$, which means that $\int {\mathcal{W}}\lnot \#\langle
f^{-1}(p)\cup \{B\}\rangle$. 

Case 2: Otherwise without loss of generality $f(W_{i})\in p$ for each $%
i<\omega $, since we are not in case 1. Define sets $V_{1}=W_{1}$, $%
V_{i}=W_{j}\cap f^{-1}(\bigcap_{k<i}f(W_{k}))$, and note that $%
\bigcup_{i<\omega }V_{i}\in \langle \mathcal{A}\cup f^{-1}(p)\rangle $ since
we are not in case 1. Then, $(f(V_{i}))_{i<\omega }$ is a decreasing
sequence, and because $f$ is finite-to-one, $(f(V_{i})\setminus
f(V_{i+1}))_{i<\omega }$ is a partition of $f(V_{1})\in p$. As $p$ is a
P-point, there exists $P\in p$ such that $P_{i}=(f(V_{i})\setminus
f(V_{i+1}))\cap P$ is finite for each $i<\omega $. Let $g:\omega \rightarrow
\omega $ be defined by $g(i)=E\left( \frac{i+1}{2}\right) $, where $E\left(
x\right) $ stands for the integer part of $x$. Let%
\begin{equation*}
R=\bigcup\nolimits_{i<\omega }(f^{-1}(P_{i})\cap \bigcup\nolimits_{j\in
\lbrack g(i),...,i]}V_{i}).
\end{equation*}%
Note that $R\cap V_{i}$ is finite for each $i<\omega $, and that%
\begin{equation*}
\limsup\nolimits_{n\in \tilde{P}}\mathrm{card\,}(f^{-1}(n)\cap R)=\infty
\end{equation*}%
for all $\tilde{P}\in p$ and $\tilde{P}\subset P$. Thus, $R^{c}\not\in
\langle \mathcal{A}\cup f^{-1}(p)\cup \widetilde{\mathcal{W}}\rangle $,
although $R^{c}\in \int \mathcal{W}$, which completes case 2. 

To complete the proof of the theorem use Corollary \ref{DDD111}.
\end{proof}

The following two, probably known, facts will be needed for Theorem %
\ref{IIIZZZ} that extends, under CH, Theorem \ref{III}.

\begin{fact}
\label{IIIAAA} Let $\mathcal{A}$ be a free family of subsets of $\omega$
such that $\mathcal{A} \cup \{F\}$ is not an ultrafilter subbase for any $F
\ddagger \mathcal{A}$. Let $\mathcal{F}$ be a \linebreak countable family of
sets such that $\mathcal{A} \ddagger \mathcal{F}$. Then $\mathcal{A} \cup
\mathcal{F}$ is not an ultrafilter subbase.
\end{fact}

\begin{proof}
Without loss of generality, we may assume, that $(F_{n})_{n<\omega }$ is a
decreasing sequence of sets, such that $F_{n+1}\not\in \langle \mathcal{A}%
\cup \{F_{n}\}\rangle $. Put $B_{n}=F_{n}\setminus F_{n+1}$ and define $%
B^1=\bigcup_{n<\omega }B_{2n}$, $B^2=\bigcup_{n<\omega }B_{2n+1}$. Clearly
at least one of sets $B^1$, $B^2$ interact $\mathcal{A}$ - say $B_1$ does.
If $B^1\not\in \langle \mathcal{A}\cup \mathcal{F}\rangle $ then we are
done. Suppose that $B^1\in \langle \mathcal{A}\cup \mathcal{F}\rangle $, and
denote by $n_{0}$, the minimal $n<\omega $ that $B^1\in \langle \mathcal{A}%
\cup \{F_{n}\}\rangle $. But $F_{n_{0}+1}\cap B^1=F_{n_{0}+2}\cap B^1$ and
so $F_{n_{0}+2}\in \langle \mathcal{A}\cup \mathcal{F}_{n_{0}+1}\rangle $,
which is a contradiction.
\end{proof}

\begin{fact}
\label{IIIBBB} Let $\mathcal{Y},\mathcal{Z}$ be free families of subsets of $%
\omega $, which are not ultrafilter subbases. If $h\in {^{\omega }\omega }$,
then there exist sets $Y$ and $Z$ such that $Y\ddagger\mathcal{Y}$, $%
Z\ddagger \mathcal{Z}$ and $h(Z)\lnot \ddagger Y$.
\end{fact}

\begin{proof}
Take any $O$ such that $O\ddagger \mathcal{Y}$ and $O^c \ddagger \mathcal{Y}$%
.

If $h^{-1}(O)\neg \ddagger \mathcal{Z}$ then $Y=O$, $Z=(h^{-1}(O))^c$;

if $h^{-1}(O^c)\neg \ddagger \mathcal{Z}$ then $Y=O^c$, $Z=h^{-1}(O)$;

if $h^{-1}(O) \ddagger \mathcal{Z}$ and $h^{-1}(O^c) \ddagger \mathcal{Z}$
then $Y=O$, $Z=(h^{-1}(O))^c$.
\end{proof}

\begin{thm}
\label{IIIZZZ} (CH) 
If $p$ is a P-point, then there exists a set $\mathfrak{U}$ of cardinality $%
\mathfrak{c}$ of Rudin-Keisler incomparable P-points with $u >_{RK} p$ for
each $u\in \mathfrak{U}$.
\end{thm}

\begin{proof}
First repeat the proof of Theorem \ref{III} except for the last line. Then
continue as follows.

We range all contours in a sequence $(\int {\mathcal{W}}_{\alpha })_{\alpha <%
\mathfrak{b}}$ and ${^{\omega }\omega }$ in a sequence $(f_{\beta })_{\beta <%
\mathfrak{b}}$. We will build a family $\{({\mathcal{F}}_{\alpha }^{\beta
})_{\alpha <\mathfrak{b}}\}_{\beta <\mathfrak{b}}$ of increasing $\mathfrak{b%
}$-sequences $({\mathcal{F}}_{\alpha }^{\beta })_{\alpha <\mathfrak{b}}$ of
filters such that:

1) Each ${\mathcal{F}}_{\alpha }^{\beta }$ is generated by $\mathcal{A}$
together with some family of sets of cardinality $<\mathfrak{b}$;

2) ${\mathcal{F}}^\beta_0 = \mathcal{A}$ for each $\beta<\mathfrak{b}$;

3) For each $\alpha ,\beta <\mathfrak{b}$, there exists $F\in {\mathcal{F}}%
_{\alpha +1}^{\beta }$ such that $F^{c}\in \int {\mathcal{W}}_{\alpha }$;

4) For every limit $\alpha $ and for each $\beta $, let ${\mathcal{F}}%
_{\alpha }^{\beta }=\bigcup_{\gamma <\alpha }{\mathcal{F}}_{\gamma }^{\beta
} $;%

5) For each $\alpha ,\gamma <\alpha ,\beta _{1}$, $\beta _{2}<\alpha ,$
there exists a set $F\in {\mathcal{F}}_{\alpha +1}^{\beta _{1}}$ such that $%
(f_{\gamma }(F))^{c}\in {\mathcal{F}}_{\alpha +1}^{\beta _{2}}$.

The existence of such families follows by a standard induction with
sub-inductions using Theorem \ref{DDD} and Fact \ref{IIIBBB} for Condition
5. It follows from the proof of Fact \ref{IIIAAA} that ${\mathcal{F}}%
_{\alpha }^{\beta }$ is not an ultrafilter-base for each every $\alpha $ and
$\beta $. It suffices now to take for each $\beta <\mathfrak{c}$, any
ultrafilter extending $\bigcup_{\beta <\mathfrak{c}}{\mathcal{F}}_{\alpha
}^{\beta }$ and note that, by Proposition \ref{A0}, it is a P-point.
\end{proof}

A. Blass \cite[Theorem 7]{Blass1} also proved that, under $\mathfrak{p}=%
\mathfrak{c}$, each RK-increasing sequence of P-points is upper bounded by a
P-point. By $\Level_n(T)$ we denote level $n$ in the tree $T$.

\begin{thm}
\label{JJJ} ($\mathfrak{b}=\mathfrak{c}$) If $(p_{n})_{n<\omega }$ is an
RK-increasing sequence of P-points, then there exists a P-point $u$ such
that $u>_{RK}p_{n}$ for each $n<\omega $.
\end{thm}

\begin{proof}
For each $n<\omega$ we let $f_{n}$ to be a finite-to-one function that witnesses $%
p_{n+1}>_{RK}p_{n}$. Consider a sequence $(T_n)_{n<\omega}$ of disjoint
trees such that for each $n<\omega$

1) the root of $T_n$ is equal to $n\in \omega$;

2) $\mathrm{Level\,}_m(T_n) = \{f_m^{-1}(k):$ $k\in \mathrm{Level\,}_{m-1}(T_n)\}$ for $1\leq m\leq n$;

3) $\mathrm{Level\,}_m(T_n)=\emptyset$  for $m>n$.

Since $L_\infty=\bigcup_{n<\omega}\mathrm{max\,} T_n$ is countably infinite
we treat it as $\omega$ as well as $L_m=\bigcup_{n<\omega}\mathrm{Level\,}%
_m(T_n)$. Let $g_m: L_\infty\setminus \bigcup_{k<m}L_k$ be a function
defined by order of the trees $T_n$.

On $L_\infty$ we define a family of sets: $\mathcal{B}=\bigcup_{n<\omega}
g_m^{-1}(p_m)$.

To conclude it suffices, by Corollary \ref{DDD111}, to show that $\mathcal{B}
$ is not quasi-finer than any contour, thus by Theorem \ref{DDD} it suffices
to prove that $g_{m}^{-1}(p_{m})$ is not quasi-finer than any contour, for
any $m$. But this is an easier version of the fact which we established in the proof of
Theorem \ref{III}.
\end{proof}

In \cite{Blass1}, A. Blass asked (Question 4) which ordinals could be
embedded in the set of P-points, noticing that such an ordinal could not be
greater than $\mathfrak{c}^+$. The question was also considered by D.
Raghavan and S. Shelah in \cite{Raghavan-Shelah} and answered, under MA, by
B. Kuzeljevic and D. Raghavan in a recent paper \cite{Kuz-Rag}.

We prove that, under $\mathfrak{b}=\mathfrak{c}$, there is an order
embedding of each ordinal less than $\mathfrak{c}^{+}$ into P-points. To
this end, we need some (probably known) facts associated with the following
definition: we say that a subset $A$ of ${^{\omega }\omega }$ is \textit{%
sparse} if $\mathrm{lim\,}_{n<\omega }\left\vert f(n)-g(n)\right\vert
=\infty $ for each $f,g\in A$ such that $f\not=g$.

\begin{fact}
\label{LLL} There exists a strictly $<^{\ast }$-increasing sparse sequence $%
\mathcal{F}=(f_{\alpha })_{\alpha <\mathfrak{b}}\subset {^{\omega }\omega }$
of nondecreasing functions such that $f_{\alpha }\left( n\right) \leq n$ for
each $n<\omega $ and $\alpha <\mathfrak{b}$.
\end{fact}

\begin{proof}
First, we build, from the definition of $\mathfrak{b}$, an $<^{\ast }$%
-increasing sparse sequence $(g_{\alpha })_{\alpha <\mathfrak{b}}\subset {%
^{\omega }\omega }$ of nondecreasing functions that fulfill the following
condition: if $\alpha <\beta <\mathfrak{b}$, then $g_{\alpha }(n)>g_{\beta
}(2n)+n$ for almost all $n<\omega $.
Then a $\mathfrak{b}$-sequence $(f_\alpha)_{\alpha<\mathfrak{b}}$ defined by $%
f_\alpha(m)=m - \mathrm{max\,}\{n: g_\alpha(n)<m\}$ is as desired.
\end{proof}

\begin{fact}
If an ordinal number $\alpha $ can be sparsely embedded in ${^{\omega
}\omega }$, under identity, as nondecreasing functions, then $\alpha $ can
be sparsely embedded in ${^{\omega }\omega }$, under any function $f$ that
converges to $\infty$.
\end{fact}

\begin{proof}
Let $h<f$ be a nondecreasing function such that $h(n+1)-h(n)\leq 1$. Let $%
(g_\beta)_{\beta<\alpha}$ be an embedding of $\alpha$. Define $%
(f_\beta)_{\beta<\alpha}$ by: $f(\alpha)(k)=g_\alpha(n)$ if and only if $h(k)=n$.
\end{proof}

\begin{fact}
\label{LLL111} If an ordinal number $\alpha $ can be sparsely embedded in ${%
^{\omega }\omega }$ as nondecreasing functions that are less than any
function $f\in {^{\omega }\omega }$, then $\alpha $ can be sparsely embedded
as nondecreasing functions between any sparse pair of functions $g<^{\ast
}h\in {^{\omega }\omega }$.
\end{fact}

\begin{proof}
Without loss of generality,  we assume $f$ to be nondecreasing. Let $(f_{\beta })_{\beta
<\alpha }$ be an embedding of $\alpha $ under $f$. Clearly, it suffices to
prove that there is an embedding under $s$ defined by $s(n)=h(n)-g(n)$ if $%
h(n)\geq g(n)$ and $s(n)=0$ otherwise.

Define a sequence $(k(n))_{n<\omega }$ by $k(0)=\mathrm{min\,}\{m:s(i)\geq
f(0)$ for all $i\geq m\}$, $k(n+1)=\mathrm{min\,}\{m:m>k(n)\text{ }\&\text{ }%
s(i)\geq f(n+1)$ for all $i\geq m\}$.
Finally define $g_\alpha$ as follows: $g_\alpha(n)=f_\alpha(m)$ if and only if $k(n) \leq m <
k(n+1)$.
\end{proof}

\begin{fact}
\label{LLL444} For each $\gamma <\mathfrak{b^{+}}$, there exists a strictly $%
<^{\ast }$-increasing sparse sequence $\mathcal{F}=(f_{\alpha })_{\alpha
<\gamma }\subset {^{\omega }\omega }$ of nondecreasing functions.
\end{fact}

\begin{proof}
Facts \ref{LLL} and \ref{LLL111} clearly imply that the first ordinal number
which cannot be embedded as a sparse sequence in ${^{\omega }\omega } $
under $\mathrm{id\,}_{\omega }$ is equal to $\alpha$ or to $\alpha +1$ where $\alpha$ is a limit number.
Facts \ref{LLL} and \ref{LLL111} also imply that the set of ordinals
less than $\alpha$ is closed under $%
\mathfrak{b}$ sums.

Indeed, let $\beta$ be the minimal ordinal number $< \mathfrak{b}^+$ that may not be embedded
under identity as an $<^*$-increasing sparse sequence. Clearly $\mathrm{cof\,%
}(\beta) \leq \mathfrak{b}$. Take an increasing sequence $(\alpha_\delta)_{%
\delta< \mathrm{cof\,} \beta}$ that converges to $\beta$. Clearly for each $%
\alpha< \beta$ there is $(g^\alpha_\eta)_{\eta<\alpha}$ - an embedding of $%
\alpha$ into ${^\omega \omega}$ as a sparse sequence under identity. By Fact %
\ref{LLL111} for each $\alpha< \beta$ there is an $<^*$ - increasing sparse
sequence of $(f^\alpha_\eta)_{\eta<\alpha}$ such that $f_\alpha <^*
f^\alpha_\eta <^* f_{\alpha+1}$ (for $f_\alpha$, $f_{\alpha+1}$ from the
proof of Fact \ref{LLL}). Now $(f^\alpha_\eta)_{\alpha<\mathrm{cof\,}%
(\beta), \eta<\alpha}$ with lexicographic order is a required embedding of $%
\beta$.

Thus, this number is not less than $\mathfrak{b}^{+}$.
\end{proof}

\begin{thm}
\label{M0} ($\mathfrak{b}=\mathfrak{c}$) For each $\gamma <\mathfrak{b^{+}}$%
, for each P-point $p$ there exists an RK-increasing sequence $\{p_{\alpha
}:\alpha <\gamma\}$ of P-points such that $p_{0}=p$.
\end{thm}

\begin{proof}
Note that $\mathrm{cof\,}(\gamma )\leq \mathfrak{b}$. Consider a set of
pairwise disjoint trees $T_{n}$ such that each $T_{n}$ has a minimal
element, each element of $T_{n}$ has exactly $n$ immediate successors, and
each branch has the highest $\omega $.

Let $\{f_{\alpha }\}_{\alpha <\gamma}\subset {^{\omega }\omega }$ be a
sparse, strictly $<^{\ast }$-increasing sequence, the existence of which is
demonstrated by Fact \ref{LLL444}. For each $\alpha <\gamma $, define%
\begin{equation*}
X_{\alpha }=\bigcup\nolimits_{n<\omega }\mathrm{Level\,}_{f_{\alpha
}(n)}T_{n}.
\end{equation*}%
For each $\alpha <\beta \leq \gamma$, define%
\begin{equation*}
f_{\alpha }^{\beta }:\bigcup\nolimits_{\{n<\omega: f_{\alpha }(n)<f_{\beta
}(n)\}}\mathrm{Level\,}_{f_{\beta }(n)}T_{n}\rightarrow
\bigcup\nolimits_{\{n<\omega: f_{\alpha }(n)<f_{\beta }(n)\}}\mathrm{Level\,}%
_{f_{\alpha }(n)}T_{n}
\end{equation*}%
that agrees with the order of trees $T_{n}$ for $n<\omega $ such that $%
f_{\alpha }(n)<f_{\beta }(n)$. Note that $\mathrm{dom\,}f_{\alpha }^{\beta }$
is cofinite on $X_{\beta }$ for each $\alpha <\beta $.

Let $p=p_{0}$ be a P-point on $X_{0}=\bigcup_{n<\omega }\mathrm{Level\,}%
_{0}T_{n}$. We proceed by recursively building a filter $p_{\beta }$ on $%
X_{\beta }$. Suppose that $p_{\alpha }$ are already defined for $\alpha
<\beta $. If $\beta$ is a successor, then it suffice to repeat a proof of Theorem \ref%
{III} for $P_{\beta-1}$ and $f^{\beta}_{\beta-1}$.

So suppose that $\beta$ is limit. Let $R\subset \beta $ be cofinite with $%
\beta$ and of order type less than or equal to $\mathfrak{b}$. Define a
family
\begin{equation*}
\mathcal{C}=\bigcup_{\alpha \in R }\{(f_{\alpha }^{\beta })^{-1}(p_{\alpha
}) \},
\end{equation*}%
which is obviously free.

Clearly each filter that extends $\mathcal{C}$ is $RK$-greater than each $%
p_{\alpha }$ for $\alpha <\beta $. But we need a P-point extension. Thus, by
Corollary \ref{DDD111} it suffices to prove that $\mathcal{C}$ is not
quasi-finer than any contour. Thus, by Theorem \ref{DDD} it suffices to prove
that $\bigcup_{\gamma \in R,\gamma \leq \alpha }\{(f_{\gamma }^{\beta
})^{-1}(p_{\gamma })\}\subset (f_{\alpha }^{\beta })^{-1}(p_{\alpha })$ is
not quasi-finer than any contour, for each $\alpha \in R$. But it is (an
easier version of) what we did in the proof of Theorem \ref{III}.
\end{proof}

%
%
%

By Theorem \ref{JJJ} and \ref{M0} the following natural question arises:

\begin{quest}
\label{M1} What is the least ordinal $\alpha $ such that there exists an
unbounded embedding of $\alpha $ into the set of P-points  \footnote{
A recent paper by D. Raghavan and J. L. Verner \cite{Rag-Ver}
showed that under $\lozenge$, the cardinal $\omega_1$ is the answer, but we
steal do not know the answer in terms of cardinal invariants which, as we
suppose, play an important role in this domain.}?
\end{quest}

A. Blass \cite[Theorem 8]{Blass1} also proved that, under $\mathfrak{p}=%
\mathfrak{c}$, there is an order-embedding of the real line into the set of
P-points. We will prove the same fact, but under $\mathfrak{b}=\mathfrak{c}$%
. Our proof is based on the original idea of set X defined by A. Blass.
Therefore, we quote the beginning of his proof, and then use our machinery.

\begin{thm}
\label{NNN} ($\mathfrak{b}=\mathfrak{c}$) There exists an order-embedding of
the real line into the set of P-points.
\end{thm}

\begin{proof}
------------------------ (beginning of quotation) -------------------------

Let $X$ be a set of all functions $x: \mathbb{Q} \to \omega$ such that $%
x(r)=0$ for all but finitely many $r\in \mathbb{Q}$; here $\mathbb{Q}$ is
the set of rational numbers. As X is denumerable, we may identify it with $%
\omega$ via some bijection. For each $\xi \in \mathbb{R}$, we define $h_\xi:
X \to X$ by

\begin{equation*}
h_{\xi }(x)(r)=\left\{
\begin{array}{c}
x(r)\text{ if }r<\xi , \\
0\text{ if }r\geq \xi .%
\end{array}%
\right.
\end{equation*}

Clearly, if $\xi < \eta$, then $h_\xi \circ h_\eta = h_\eta \circ h_\xi =
h_\xi$. The embedding of $\mathbb{R}$ into P-points will be defined by $\xi
\to D_\xi = h_\xi(D)$ for a certain ultrafilter $D$ on $X$. If $\xi < \eta$,
then

$D_{\xi }=h_{\xi }(D)=h_{\xi }\circ h_{\eta }(D)=h_{\xi }(D_{\eta })\leq
D_{\eta }$.

We wish to choose $D$ in such a way that

(a) $D_\xi \not\cong D_\eta$ (therefore, $D_\xi < D_\eta$ when $\xi < \eta$%
), and

(b) $D_\xi$ is a P-point.

Observe that it will be sufficient to choose $D$ such that

(a') $D_\xi \not\cong D_\eta$ when $\xi < \eta$ and both $\xi$ and $\eta$
are rational, and

(b') $D$ is a P-point.

Indeed, (a') implies (a) because $\mathbb{Q}$ is dense in $\mathbb{R}$. If
(a) holds, then $D_{\xi-1} < D_\xi$, so $D_\xi$ is a nonprincipal
ultrafilter $\leq D$; hence (b') implies condition (b).

Condition (a') means that, for all $\xi <\eta \in \mathbb{Q}$ and all $%
g:X\rightarrow X$, $D_{\eta }\not=g(D_{\xi })=gh_{\xi }(D_{\eta })$. By
Theorem \ref{MMM111}, this is equivalent to $\{x:gh_{\xi }(x)=x\}\not\in
D_{\eta }$, or by our definition of $D_{\eta }$,
\begin{equation}
\{x:gh_{\xi }(x)=h_{\eta }(x)\}=h_{\eta }^{-1}\{x:gh_{\xi }(x)=x\}\not\in D.
\tag{a\textquotedblright }  \label{a}
\end{equation}

We now proceed to construct a P-point $D$ satisfying (\ref{a}) for all $\xi
<\eta \in \mathbb{Q}$ and for all $g:X\rightarrow X$; this will suffice to
establish the theorem.

---------------------------- (end of quotation)
-------------------------------

List all pairs $(\xi, \eta)$, $\xi < \eta \in \mathbb{Q}$ in the sequence $%
(\xi_i, \eta_i)_{i<\omega}$. For each $g\in {^{X}X}$, $\xi <\eta \in \mathbb{%
Q}$, define ${A}_{g,\xi ,\eta }=\{x\in X:gh_{\xi }(x)\not=h_{\eta }(x)\}$, $%
\mathcal{A}_{i}=\mathcal{A}_{\xi _{i},\eta _{i}}=\{{A}_{g,\xi _{i},\eta
_{i}}:g\in {^{X}X}\}$, and $\mathcal{A}=\bigcup_{i<\omega }\mathcal{A}_{i}$.
Clearly, $\mathcal{A}$ is free. \vspace{0mm}

We claim that $\mathcal{A}$ is not quasi-finer than any contour.

Indeed, by Theorem \ref{DDD}, it suffices to prove that for each $n<\omega $%
, $\bigcup_{i<n}\mathcal{A}_{i}$ is not quasi-finer than any contour.
Suppose not and take (by Remark \ref{AAA111}) the witnesses $i_{0}$ and $%
\int \mathcal{W}$ such that $\int {\mathcal{W}}\subset \langle
\bigcup_{i<i_{0}}\mathcal{A}_{i}\cup \widetilde{\mathcal{W}}\rangle $. For
each $n<\omega $, consider the condition $(S_{n})$:
\begin{equation*}
\exists x_{n}:\forall i<i_{0}:\left( x_{n}\in h_{\xi _{i}}({\widetilde{W}_{1}%
})\text{ }\&\text{ }\mathrm{card\,}(h_{\eta _{i}}(h_{\xi
_{i}}^{-1}(x_{n}))\cap h_{\eta _{i}}({\widetilde{W}_{n}}))>n\right) .
\end{equation*}

Case 1: $S_{n}$ is fulfilled for all $n<\omega $. Then, for each $n<\omega $%
, $j<n$, choose $x_{n}^{j}\in \widetilde{W}_{n}$ such that $h_{\xi_i
}(x_{n}^{j})=x_{n}$ and $h_{\eta _{i}}(x_{n}^{j_{0}})\not=h_{\eta
_{i}}(x_{n}^{j_{1}})$ for $j_{0}\not=j_{1}$. Define ${E}=\bigcup_{n<\omega
}\bigcup_{j\leq n}\{x_{n}^{j}\}$. Clearly ${E}^{c}\in \int {\mathcal{W}}$,
but ${E}\not\subset \bigcup_{i<i_{0}}\bigcup_{g\in \mathcal{G}}({A}_{g,\xi
_{i},\eta _{i}})\cup \bigcup_{l<m}W_{l}$ for any choice of finite family $%
\mathcal{G}\subset {^{X}X}$ and for any $m<\omega $.

Case 2: $S_{n}$ is not fulfilled for some $n_{0}<\omega $. Then, there exist
functions $\{g_{n,i}\}_{n\leq n_{0},i<i_{0}}\subset {^{X}X}$ such that $%
\widetilde{W}_{1}\subset \bigcup_{n\leq n_{0}}\bigcup_{i<i_{0}}({A}%
_{g_{n,i},\xi _{i},\eta _{i}})\cup \bigcup_{n\leq n_{0}}W_{n}$, i.e., $\int
\mathcal{W}\lnot \#\langle \bigcup_{i<i_{0}}\mathcal{A}_{i}\cup \widetilde{%
\mathcal{W}}\rangle $. 
\vspace{2mm}

%

Corollary \ref{DDD111} completes the proof.
\end{proof}


The \textit{long line} is defined as $\mathbb{L} = \omega_1 \times (0,1]$
ordered lexicographically. If $f:Y \to \omega$, then the \textit{support} of
$f$ is defined as
\begin{equation*}
\mathrm{supp\,}(f)=\{y\in Y: f(y) \not=0\}.
\end{equation*}

\begin{lemm}
\label{OOO} ($\mathfrak{b}=\mathfrak{c}$) For each P-point $p$, there exists
an order-embedding of the real line into the set of P-points above $p$.
\end{lemm}

\begin{proof}
We will combine ideas form proofs of Theorems \ref{III} and \ref{NNN} with
some new arguments.

Again, let $X$ be a set of all functions $x: \mathbb{Q} \to \omega$ such
that $x(r)=0$ for all but finitely many $r\in \mathbb{Q}$. Since $X$ is
infinitely countable, we treat it as $\omega$.

Let $p$ be a P-point on $X$ such that for each $q\in \mathbb{Q}$ and for each $%
P\in p$ there exists $x \in P$ such that $\mathrm{max\,} \mathrm{supp\,}(x)<
q$. Let $f\in {^{X}X }$ be a finite-to-one function such that $%
\limsup\nolimits_{x\in P}\mathrm{card\,}(f^{-1}(x))=\infty$ for all $P\in p$
and that $\mathrm{max\,} \mathrm{supp\,} x < \mathrm{max\,} \mathrm{supp\,}
f(x)$. Again, we define a family $\mathcal{A}$ as follows: $A\in \mathcal{A}$
if and only if there exist $i<\omega $ and $P\in p$ such that $\mathrm{card\,}%
(f^{-1}(x)\setminus A)\leq i$ for each $x\in P$. For each $\xi \in \mathbb{R}
$, we again define functions $h_\xi: X \to X$ by

\begin{equation*}
h_{\xi }(x)(r)=\left\{
\begin{array}{c}
x(r)\text{ if }r<\xi , \\
0\text{ if }r\geq \xi .%
\end{array}%
\right.
\end{equation*}

List all rational numbers in $\omega$-sequence $\mathbb{Q}=(q_i)_{i<\omega}$%
. Let $\mathcal{B}_i=h^{-1}_{q_i}(\mathcal{A} \cup f^{-1}(p))$ and let $%
\mathcal{B} = \bigcup_{i<\omega} \mathcal{B}_i$.

Our aim is to prove that $\mathcal{B}$ can be extended to such a P-point $Q$
that $h_\xi(Q) \not= h_\eta(Q)$ for each $\xi\not= \eta \in \mathbb{Q}$ (and
thus for each $\xi\not= \eta \in \mathbb{R}$).\nolinebreak\hspace{19mm}%
\nolinebreak(4)

To this end, we add to $\mathcal{B}$ a family $\mathcal{C}$ defined as
follows: list all pairs $(\xi, \eta)$, $\xi < \eta \in \mathbb{Q}$ in the
sequence $(\xi_i, \eta_i)_{i<\omega}$. For each $g\in {^{X}X}$, $\xi <\eta
\in \mathbb{Q}$, define ${C}_{g,\xi ,\eta }=\{x\in X:gh_{\xi
}(x)\not=h_{\eta }(x)\}$, $\mathcal{C}_{i}=\mathcal{C}_{\xi _{i},\eta
_{i}}=\{{C}_{g,\xi _{i},\eta _{i}}:g\in {^{X}X}\}$, and $\mathcal{C}%
=\bigcup_{i<\omega }\mathcal{C}_{i}$.

Thus to prove (4), it suffices by Corollary \ref{DDD111} to prove that $%
\mathcal{B}\cup \mathcal{C}$ is not quasi-finer than any contour. Thus, by
Theorem \ref{DDD}, in order to prove (4), it suffices to prove that: 

$\mathcal{D}_{i}$ is not quasi-finer than any contour for $\mathcal{D}_{i}$
defined for $i<\omega $ as follows:
\begin{equation*}
\mathcal{D}_{i}=\bigcup_{j\leq i}\mathcal{B}_{i}\cup \bigcup_{j\leq i}%
\mathcal{C}_{i}.
\end{equation*}

\vspace{-12mm}\hspace{127mm}(5)\vspace{10mm}

First, to prove it, we notice that $\mathcal{D}_i$ is free. Indeed, define $%
q_m = \mathrm{min\,} \{q_j: j \leq i\}$, $q_M = \mathrm{max\,} \{q_j: j \leq
i\}$ and $\xi_m=\mathrm{min\,}\{\xi_j: j\leq i\} $ and note that $%
h^{-1}_{q_M}(x) \subset h^{-1}_{q_j}(x)$ for each $j \leq i$ and for each $x$
such that $\mathrm{max\,}\mathrm{supp\,}(x) < \mathrm{min\,}\{ q_m, \xi_m \}$%
. It is easy to see that $h^{-1}_{q_M}(x) \cup \bigcup_{j\leq i}\mathcal{C}%
_i $ is free, hence $\mathcal{D}_i$ is free.

\vspace{1mm} Fix $i$.
and suppose that (5) does not hold. So take a witness $\int \mathcal{W}$.
From Remark \ref{AAA111}, without loss of generality, we may assume that $%
\int {\mathcal{W}}\subset \langle \mathcal{D}_{i}\cup \widetilde{\mathcal{W}}%
\rangle $.

Let $A_W \in \mathcal{A}$, $P_W\in p$, $n_W\in \omega$, $C_n^W \in \mathcal{C%
}$ for $n\leq n_W$ and $l_W \in \omega$. Define
\begin{equation*}
W^*(A_W, P_W, C_1^W, \ldots , C_{n_W}^W, \widetilde{W_{l_W}}) =
\bigcap_{n\leq n_W} C_n^W \cap \bigcap_{j \leq i} h^{-1}_{\xi_j} (A_W \cap
P_W) \cap \widetilde{W_{l_W}} .
\end{equation*}

Define $W \in \mathcal{W}^-$ if and only if $W \in \int \mathcal{W}$ and $W$ is
co-finite or empty on each $W_i$. We will say that a set $W\in \mathcal{W}$
is $\mathit{attainable}$ (by $(n_A, P_W, n_W, l_W)$) if there exist $A_W \in
\mathcal{A}_{n_A}$, $\{C_k^W \in \bigcup_{j \leq i} \mathcal{C}_j: k \leq
n_W\}$ such that 
the condition $W^*(A_W, P_W, C_1^W, \ldots , C_{n_W}^W, \widetilde{W_{l_W}%
}) \subset W$ is satisfied. The complement (to $\widetilde{W_1}$) of the attainable set
is called \textit{removable} and sometimes we point out by which variables,
sets, functions.

Since $\int \mathcal{W} \subset \langle \mathcal{D}_i \rangle$, thus
each set $W \in \mathcal{W}^-$ is attainable.

Consider a sequence of possibilities:

1) $l$ cannot be fixed, i.e., for each $l\in\omega$ there exists $W \in {%
\mathcal{W}^-}$ such that $W$ is not removable by any $(n_A, P_W, n_W, l)$;

2) $l$ can be fixed, but $n_A$ cannot, i.e., for each $W\in \mathcal{W}^-$,
$W$ is attainable by some $(n^A_W, P_W, n^C_W, l)$, but for each $n_A$ there
exists $W^{\prime }\in {\mathcal{W}^-}$ such that $W^{\prime }$ is not
attainable by any $(n^A, P_{W^{\prime }}, n^C_{W^{\prime }}, l)$;

3) $l$ and $n^A$ can be fixed, but $n_C$ cannot;

4) $l$, $n^A$ and $n^C$ can be fixed, but $P$ cannot;

5) $l$, $n^A$, $n^C$ and $P$ can be fixed.

Note that each set $W \in \mathcal{W}^-$ is attainable if and only if an alternative of cases 1) to 5)
holds.

In case 1) for each $l$, let $W_l\subset \widetilde{ {W}_l}$ and $W_l \in
\mathcal{W}^-$ be a witness that $l$ may not be fixed. Note that
$\bigcup_{l<\omega}W_l \in \mathcal{W}^-$ and that $\bigcup_{l<\omega}W_l$
may not be removed by any $(n^A_W, P_W, n^C_W, l_W)$.

In case 2) we proceed like in case 1). Note that if $l^{\prime }\geq l$ and ${%
n^A}^{\prime}$ and a set $W \in \mathcal{W}^-$ is not removable by any $({%
n^A}^{\prime }, P_W, n^C_W, l^{\prime })$ then the set $W$ is also not
removable by any $(n^A, P_W, n^C_W, l)$. Thus it suffices to consider cases
when $l=n^A$. For each $l$, let $W_l\subset \widetilde{ {W}_l}$ and $W_l \in
\mathcal{W}^-$ be a witness that $l$ and $n^A=l$ may not be fixed. Again
note that 
$\bigcup_{l<\omega}W_l \in \mathcal{W}^-$ and that $\bigcup_{l<\omega}W_l$
may not be removed by any $(n^A_W, P_W, n^C_W, l_W)$.

In case 3) we proceed just like in case 2), not using that $n^A$ is fixed.

In case 4) for $k<\omega$, let $X_k$ be the set of those $x\in X$ that for
all $U\subset f^{-1}(x)$ such that $\mathrm{card\,}(f^{-1}(x) \setminus U)
\leq n^A$, for all partitions of a set $\bigcup_{j=1}^i(\xi_j^{-1}(X) \cap
\widetilde{(W_k)})$ on the sets $X_{m,n}$, for $m<n\leq i$, there exist $%
m_0, n_0$ such that $m_0<n_0\leq i$ and there exist $x_1, \ldots, x_{n^C+1}
\in X_{m_0,n_0}$ that for $\xi_{(min)}=\mathrm{min\,}\{\xi_{m_0},
\xi_{n_0}\} $, $\xi_{(max)}=\mathrm{max\,}\{\xi_{m_0}, \xi_{n_0}\}$ there is
$\xi_{(min)}(x_r)=\xi_{(min)}(x_j)$ for $r,j \leq n^C+1$ and $%
\xi_{(max)}(x_r)\not=\xi_{(max)}(x_j)$ for $r,j \leq n^C+1$, $r\not= j$.

Clearly, $(X_{k})$ is a decreasing sequence. If there exists $k$ such that $%
X_{k}\not\in p$ then putting $l=k$ there exists a set $P=(X_{k})^{c}$ such
that all $W\in \mathcal{W}^{-}$ may be attained by $(n^{A},P,n^{C},l)$ so we
would be in case 6, so, without loss of generality, $X_{k}\in P$ for each $%
k<\omega $.

Thus take a partition of $X$ by $(X_k\setminus X_{k+1})$. Since $p$ is a
P-point, and since $X_k\in p$ thus there exists $P_0\in p$ such that $%
P_k=P_0\cap (X_k \setminus X_{k+1})$ is finite for all $k<\omega$. For each $%
x\in P_k$ there exists a finite(!) set $K_{k,x}\subset \widetilde{(W_k)}
\cap \bigcup_{j\leq i}h^{-1}_{\xi_i}$ that may not be removed by $%
(n^A,X,n^C,k)$. (The proof that $K_{k,x}$ may be chosen finite is
analogical, but easier, to that of case 5)). Take $%
K=\bigcup_{k<\omega, x\in P_k}K_{k,x}$ and notice that $\widetilde{(W_1)}%
\setminus K \in \mathcal{W}^-$ and that $\widetilde{(W_1)}\setminus K$ is
not removable by $(n^A, P, n^C, l)$ for any $P\in p$.

In case 5) arrange $\bigcup_{j \leq i} h^{-1}_{\xi_j} (P) \cap \widetilde{%
W_{l}}$ into a sequence $(x_k)_{k<\omega}$. Let $R(x) = {\binom{{\mathrm{%
card\,}(f^{-1}(x))} }{{\mathrm{card\,}(f^{-1}(x)) - n^C}}}$, where $\binom{{n%
} }{{k}}$ denotes a binomial coefficient, and let $(A_{x,r})_{r \leq R}$ be
a sequence of all subsets of $f^{-1}(x)$ of cardinality equal to $\mathrm{%
card\,}(f^{-1}(x) - n^C)$. Consider a tree $T$, where the root is $\emptyset$
and on a level $k$ the nodes are pairs of natural numbers $j,r$ such that $%
j\leq i$ and $r\leq R(f(x_k))$ and, for each branch $\tilde{T}$ of $T$, $%
\pi_2(\tilde{T}(k_1))=\pi_2(\tilde{T}(k_2))$ if $f(x_{k_1})=f(x_{k_2})$,
where $\tilde{T}(k)$ is an element of level $k$ of a branch $\tilde{T}$ and $%
\pi_2$ is a projection on the second coordinate. We see $j$ as a choice to
which class $\mathcal{C}_j$ does a set $C_{(.)}^{(.)}$ belongs and we see $r$
as a choice of one of sets $A_{f(x),r}$ that $C_{(.)}^{(.)}$ together with $%
A_{f(x),r}$ removes $x_k$.

Clearly, the complements of all finite sets belong to $\int \mathcal{W}$, so
each finite set is removable. \hspace{103mm} (6)

The maximal element of the branch $\tilde{T}$ has no successors if and only if there is $j
\leq i$ such that there is no $n$ sets in $\mathcal{C}_j$ that remove all $%
x_t$ such that $\tilde{T}(k)=j$ and $f(x_k)\in A_{f(x_k), \pi_2(\tilde{T}%
(k))}$. It implies that the set $\{x_k: \tilde{T}(k)=j,$ $f(x_k)\in
A_{f(x_k), \pi_2(\tilde{T}(k))}\}$ contains more than $n$ different
elements, say $x^1, \ldots x^{n+1}$, such that $h_{\xi_j}(x_{s_1}) =
h_{\xi_j}(x_{s_2})$ and $h_{\eta_j}(x_{s_1}) \not= h_{\eta_j}(x_{s_2})$ for $%
s_1\not= s_2$, $s_1, s_2 \in \{1, \ldots, n+1\}$.

By K$\ddot{\mathrm{o}}$nig Lemma if all branches are finite, then the height
of the tree $T$ is finite, and so there are irremovable finite sets in
contrary to (6). Thus there is infinite branch and the whole set $\bigcap_{j
\leq i} h^{-1}_{\xi_i} (P) \cap \widetilde{W_{l}}$ is removable.
\end{proof}


%
%
%
%

\vspace{5mm}

As an immediate consequence of Lemma \ref{OOO} (with the use of Theorem \ref%
{JJJ}) we have the following:

\begin{thm}
\label{PPP} ($\mathfrak{b}=\mathfrak{c}$) For each P-point $p$, there exists
an order embedding of the long line into the set of P-points above $p$.
\end{thm}

\begin{rem}
Note that there is a potential chance to improve Theorem \ref{PPP} in the
virtue of Question \ref{M1}, i.e. if, in some model, for each $\alpha <
\kappa$ (for some cardinal invariant $\kappa$) each RK-increasing $\alpha$%
-sequence of P-points is upper bounded by a P-point, then (in that model) if
$\mathfrak{b}=\mathfrak{c}$, then, above each P-point, there is an order
embedding of a $\kappa$-long-line into the set of P-points.
\end{rem}

\section{Cardinal $\mathfrak{q}$}

An inspection of our proofs indicates a possibility of refinement of most
results with the aid of an, {\it a priori}, new cardinal invariant.
We define $\mathfrak{q}$ to be the minimal cardinality of families $\mathcal{%
B}$, for which there exists a family $\mathcal{A}$ such that $\langle
\mathcal{A}\cup \mathcal{B}\rangle $ includes a contour, and $\langle
\mathcal{A}\cup \mathcal{C}\rangle $ includes no contour for every countable
family $\mathcal{C}$.

If $\mathbb{P}$ is a collection of families such that $\mathcal{P}\in
\mathbb{P}$ whenever $\langle \mathcal{P}\rangle $ includes a contour, then $%
\mathfrak{q}$ fulfills

\begin{equation*}
\mathfrak{q}=\mathrm{min\,}\{\mathrm{card\,}(\mathcal{B}):\underset{\mathcal{%
A}}{\exists }\;\mathcal{A}\cup \mathcal{B}\in \mathbb{P}\wedge \underset{%
\mathcal{C}}{\forall }\left( \mathrm{card\,}(\mathcal{C})\leq \aleph
_{0}\Longrightarrow \mathcal{A}\cup \mathcal{C}\notin \mathbb{P}\right) \}.
\label{defq}
\end{equation*}

{Each contour has a base of cardinality $\mathfrak{d}$ \cite[Theorem 5.2]%
{Star-Ord-V-P}, which, by the way, is the minimal cardinality of bases of
contours.} Therefore, taking into account Theorem \ref{DDD}, we have

\begin{thm}
$\mathfrak{b}\leq \cof(\mathfrak{q}) \leq \mathfrak{q}\leq \mathfrak{d}$.
\end{thm}


Using the cardinal $\mathfrak{q}$, we are in a position to formulate
stronger versions, if $\mathfrak{b}<\mathfrak{q}$ is consistent, of several
of our theorems with almost unchanged proofs. Indeed, by the proof of
Theorem \ref{III} we get the following theorem:

%

\begin{thm}
($\mathfrak{q}=\mathfrak{c}$) For each P-point $p$ there exists a P-point $q$ strictly
RK-greater than $p$.
\end{thm}

By the proof of Theorem \ref{JJJ}, we have

\begin{thm}
($\mathfrak{q}=\mathfrak{c}$) If $(p_{n})_{n<\omega }$ is an RK-increasing
sequence of P-points, then there exists a P-point $u$ such that $%
u>_{RK}p_{n} $ for each $n<\omega $.
\end{thm}

By the proof of Theorem \ref{M0}, we get

\begin{thm}
($\mathfrak{q}=\mathfrak{c}$) For each P-point $p$, for each $\alpha<%
\mathfrak{b}^+$, there exists an
order embedding of $\alpha$ into P-points above $p$. 
\end{thm}

By the proof of Theorem \ref{NNN}, we obtain

\begin{thm}
($\mathfrak{q}=\mathfrak{c}$) Above each P-point, there exists an
order-embedding of the real line in the set of P-points.
\end{thm}

By the proof of Theorem \ref{PPP}, we have

\begin{thm}
\label{QQQ} ($\mathfrak{q}=\mathfrak{c}$) Above each P-point, there exists
an order embedding of the long line into the set of P-points above $p$.
\end{thm}

A relative importance of the facts formulated above depends on answers to
the following quest.

\begin{quest}
Is $\mathfrak{q}$ equal to any already defined cardinal invariant? Is $%
\mathfrak{b}<\mathfrak{q}$ consistent? Is $\mathfrak{q}<\mathfrak{d}$
consistent?
\end{quest}

\section{Variants of invariants}

The cardinal $\mathfrak{q}$ can be seen as an instance of cardinal
invariants, which can possibly be defined in order to refine certain types
of theorems, by scrutinizing the mechanisms underlying their proofs. In our
approach, such cardinals represent \textquotedblleft
distances\textquotedblright\ between certain classes of objects. They carry
some obvious questions about their relation to the usual cardinal
invariants, and in particular to those that they are supposed to replace in
potentially refined arguments.

Let $\mathbb{S}$ and $\mathbb{T}$ be collections of families (of sets or
functions, or possibly other objects) 
such that for each $\mathcal{S}\in \mathbb{S}$ there exists $\mathcal{T}\in
\mathbb{T}$ such that $\mathcal{S}\subset \mathcal{T}$. For each $\mathcal{S}%
\in \mathbb{S}$, we define%
\begin{equation*}
\mathfrak{dist}(\mathcal{S},\mathbb{T})=\mathrm{min\,}\{\mathrm{card\,}(%
\mathcal{B}):\mathcal{S}\cup \mathcal{B}\in \mathbb{T}\}.
\end{equation*}%
Let $\mathfrak{D}(\mathbb{S},\mathbb{T})=\{\mathfrak{dist}(\mathcal{S},%
\mathbb{T}):\mathcal{S}\in \mathbb{S}\}$. As a set of cardinal numbers, $%
\mathfrak{D}(\mathbb{S},\mathbb{T})$ is well ordered, hence we can define $%
\mathfrak{dist}_{\beta }(\mathbb{S},\mathbb{T})$ to be the $\beta $-th
element of $\mathfrak{D}(\mathbb{S},\mathbb{T})$  \footnote{Equivalently, $\mathfrak{dist}_{\alpha }(\mathbb{S},%
\mathbb{T})$ can be defined recursively by $\mathfrak{dist}_{0}(\mathbb{S},%
\mathbb{T})=0$, and if $\mathfrak{dist}_{\beta }(\mathbb{S},\mathbb{T})$ is
already defined for all $\beta <\alpha $, then $\mathfrak{dist}_{\alpha }(\mathbb{S},\mathbb{T})$ is equal to
\begin{equation*}
\mathrm{min\,}\hspace{-1.3mm%
}\left\{ \mathrm{card\,}(\mathcal{B}):\underset{\mathcal{S}\in \mathbb{S}}{%
\exists }\hspace{-0.9mm}\left[ \mathcal{S}\cup \mathcal{B}\in \mathbb{P}%
\wedge \underset{\mathcal{C}}{\forall }\underset{\beta <\alpha }{\forall }%
\Big(\mathrm{card\,}(\mathcal{C})\leq \mathfrak{dist}_{\beta }(\mathbb{S},%
\mathbb{T})\Longrightarrow \mathcal{S}\cup \mathcal{C}\notin \mathbb{P}\Big)%
\right] \right\} .
\end{equation*}%
}.
Moreover, if $\alpha $ is
a limit ordinal, and the cardinals $\mathfrak{dist}_{\beta }(\mathbb{S},%
\mathbb{T})$ are defined for all $\beta <\alpha $, then $\mathfrak{dist}%
_{<\alpha }(\mathbb{S},\mathbb{T}) = \mathrm{sup\,}_{\beta <\alpha }\mathfrak{%
dist}_{\beta }(\mathbb{S},\mathbb{T})$.

In particular, if $\mathbb{S}$ denotes the collection of free families of
subsets of $\omega $, and
$\mathbb{T}$ stands for the collection of families including a subbases of a contour, then
we write $\mathfrak{q}_{\alpha }=\mathfrak{dist}_{\alpha }(\mathbb{S},%
\mathbb{T})$. In order to show that $\left( \mathfrak{q}_{\alpha }\right)
_{\alpha }$ are variants of $\mathfrak{q}$, we need the following
Alternative Theorem \cite[Theorem 3.1]{DSW1}. A relation $A\subset \left\{ \left( n,k\right) :n<\omega ,k<\omega \right\} $
is called \emph{transversal} if $A$ is infinite, and $\left\{ l:\left(
n,l\right) \in A\right\} $ and $\left\{ m:\left( m,k\right) \in A\right\} $
are at most singletons for each $n,k<\omega .$

\begin{thm}
\label{thm:alt}Let $\left( \mathcal{F}_{n}\right) _{n}$ and $\left( \mathcal{%
G}_{k}\right) _{k}$ be sequences of filters on a set $X$, and let%
\begin{equation*}
\mathcal{F}:=\bigcup\nolimits_{m<\omega }\bigcap\nolimits_{n>m}\mathcal{F}%
_{n},\text{ and }\mathcal{G}:=\bigcup\nolimits_{m<\omega
}\bigcap\nolimits_{k>m}\mathcal{G}_{k}.
\end{equation*}%
If $\mathcal{F}\#\mathcal{G}$, then the following alternative holds:  $%
\mathcal{F}_{n}\#\mathcal{G}_{k}$ for a transversal set of $\left(
n,k\right) $, or $\mathcal{F}\#\mathcal{G}_{k}$ for infinitely many $k,$ or $%
\mathcal{F}_{n}\#\mathcal{G}$ for infinitely many $n.$
\end{thm}

\begin{propo}
$\mathfrak{q}_{0}=0$, $\mathfrak{q}_{1}=1$, $\mathfrak{q}_{2}=\aleph _{0}$, $%
\mathfrak{q}_{3}=\mathfrak{q}$, and $\mathfrak{q}_\alpha \leq \mathfrak{d}$
for all $\alpha$.%
\end{propo}

\begin{proof}
By taking $\mathcal{S}\in \mathbb{T}$ and $\mathcal{B}=\emptyset ,$ we
infer that $\mathfrak{q}_{0}=0$.

To see that $\mathfrak{q}_{1}=1,$ let $A$, $B$ be disjoint countably
infinite sets. Let $\mathcal{S}_{A}$ be a contour on $A$ and let $\mathcal{S}%
_{B}$ be a cofinite filter on $B$. Define a filter $\mathcal{S}$ on $A\cup B$
by $S\in \mathcal{S}$ if and only if $S\cap A\in \mathcal{S}_{A}$ and $S\cap B\in
\mathcal{S}_{B}$. Clearly $\mathcal{S}$ is not finer than a contour (since
is RK-smaller than a cofinite filter $\mathcal{S}_{B}$), and $\mathcal{S}%
\cup \{A\}$ is a subbase of a contour.

Clearly, $\mathfrak{q}_{2}$ cannot be finite. To see that $\mathfrak{q}%
_{2}=\aleph _{0},$ let $\mathcal{W}=(W_{n})_{n<\omega }$ be a partition of $%
\omega $ into infinite sets. We define a family $\mathcal{S}$ so that $S\in
\mathcal{S}$ if and only if $S$ is cofinite on each $W_{n}$. Suppose that there is a
partition $\mathcal{V}=(V_{n})_{n<\omega }$ such that $\int {\mathcal{V}}%
\subset \langle \mathcal{S}\cup S_{0}\rangle $ for some set $S_{0}$ such
that $\mathcal{S}\cup S_{0}$ is free. Let $N_{\limfunc{fin}}=\{n<\omega :%
\mathrm{card\,}(W_{n}\cap S_{0})<\omega \}$ and $N_{\mathrm{\infty \,}%
}=\omega \setminus N_{\limfunc{fin}}$. Define $S_{\limfunc{fin} }=S_{0}\cap
\bigcup_{n\in N_{\limfunc{fin}}}W_{n}$ and $S_{\infty }=S_{0}\cap
\bigcup_{n\in N_{\infty }}W_{n}$. Since $S_{\limfunc{fin} }^{c}\in \mathcal{S%
}$, without loss of generality, we can assume that $S_{0}=S_{\infty }$ and
so without loss of generality we can assume that $S_{0}=\omega $.

Note also that $\bigcup_{n<\omega }V_{n}\cap W_{i}$ is infinite for infinitely
many $i$, and so $\int \mathcal{W}\ddagger \int \mathcal{V}$. Thus we meet
the assumptions of Theorem \ref{thm:alt}, and in each of the three cases
there exist $i,j<\omega $ such that $V_{i}\cap W_{j}$ is infinite. But $%
\omega \setminus V_{i}\in \int \mathcal{V}$ and thus $\omega \setminus
(V_{i}\cap W_{j})\in \mathcal{S}$, contrary to the definition of $\mathcal{S}
$. On the other hand, by adding $\widetilde{\mathcal{W}}$ to $\mathcal{S}$,
we obtain a subbase of $\int \mathcal{W}.$

That $\mathfrak{q}_{3}=\mathfrak{q}$, follows directly from the definition
of $\mathfrak{q}$.

Finally, $\mathfrak{q}_{\alpha }\leq \mathfrak{d}$ since each contour has a
base of cardinality $\mathfrak{d}$, as we have shown in \cite[Theorem 5.2]%
{Star-Ord-V-P}.
\end{proof}


If $\mathbb{S}$ is the collection of free families, as above, but $%
\mathbb{T}$ is the collection of free ultrafilter subbases, then ,taking $\mathcal{S}$ as an empty family, clearly
$\mathfrak{dist}(\mathcal{S}, \mathbb{T})=u$ and so
$\mathfrak{dist}_{\alpha }(\mathbb{S},\mathbb{T})=\mathfrak{u}$ for some $%
\alpha $, thus we obtain variants of $\mathfrak{u}$. By Fact \ref{IIIAAA},

\begin{fact}
$\mathfrak{u}_0=0$, $\mathfrak{u}_1=1$, $\mathfrak{u}_2 \geq \aleph_1$.
\end{fact}

By the proof of Theorem \ref{IIIZZZ}, we obtain:

\begin{thm}
($\mathfrak{q}=\mathfrak{u}_{2}=\mathfrak{c}$). If $p$ is a P-point, then
there is a set $\mathfrak{U}$ of Rudin-Keisler incomparable P-points such
that $\limfunc{card}\mathfrak{U}=\mathfrak{c}$ and $u>_{RK}p$ for each $u\in
\mathfrak{U}$.
\end{thm}

\vspace{3mm}

A similar approach can be carried out for all others cardinal invariants. Its
usefulness, however, depends on the way these cardinals are used
in specific arguments.

\vspace{3mm} {\small \textsc{\ Acknowledgements}: The author is grateful to Professor Szymon Dolecki for inspiring talks during the preparation of this paper, and to the late Professor J\'ozef Burzyk for a detailed perusing and pointing out some inaccuracies in the previous version of the paper.}


\begin{thebibliography}{10}

\bibitem{Blass1} {  A. Blass, The Rudin--Keisler ordering of
P-points. Trans. Amer. Math. Soc. 179 (1973), 145-166. }

\bibitem{Blass2} {  A. Blass, Combinatorial cardinal
characteristics of the continuum, in Handbook of Set Theory, 395-489,
Springer, Dordrecht, 2010. }

{
}

{
}

\bibitem{Dol-Myn1} {  S. Dolecki, F. Mynard, Convergence
Foundations of Topology. World Scientific Publishing, 2016. }

{
}

{\footnotesize
}

{\footnotesize
}

{\footnotesize
}

{\footnotesize
}

{\footnotesize
}

\bibitem{DSW1} {  S. Dolecki, A. Starosolski, Continuous
extension of maps between sequential cascades,  	arXiv:1901.10729}

\bibitem{vanDouwen} {  E. K. van Douwen, The integers and
topology, in Handbook of Set-Theoretic Topology, 111-167, North Holland,
1984. }

{\footnotesize
}

{\footnotesize
}

{\footnotesize
}

{\footnotesize
}

{\footnotesize
}

{\footnotesize
}

{\footnotesize
}

\bibitem{Kuz-Rag} {  B. Kuzeljevic, D. Raghavan, A long chain of
P-points, J. Math. Log. 18, 1850004 (2018). }

\bibitem{Rag-Ver} {  D. Raghavan, J. L. Verner, Chains of
P-points, arXiv:1801.02410 (2018). }

{\footnotesize
}

{\footnotesize
}



\bibitem{Raghavan-Shelah} {  D. Raghavan, S. Shelah, On
embedding certain partial orders into the P-points under Rudin--Keisler and
Tukey reducibility. Trans. Amer. Math. Soc. 369 (2017), no. 6, 4433-4455. }

{\footnotesize
}

\bibitem{MERudin} {  M. E. Rudin, Partial orders on the types of
$\beta {N}$, Trans. Amer. Math. Soc., 155 (1971), 353-362. }

{\footnotesize 
}

\bibitem{Star-P-hier} {  A. Starosolski, P-hierarchy on $%
\beta\omega$, J. Symb. Log. 73, 4 (2008), 1202-1214. }

\bibitem{Star_RK_P-hier} {  A. Starosolski, The Rudin--Keisler
ordering on P-hierarchy, arXiv:1204.4173 (2012).}

\bibitem{Star-Ord-V-P} {  {\ A. Starosolski, Ordinal
ultrafilters versus P-hierarchy, Cent. Eur. J. Math. 12 (2014), no. 1,
84-96. } }


\end{thebibliography}
\end{document}